\documentclass[a4paper,10pt]{amsart}
\usepackage{amssymb}
\usepackage{amsmath}
\usepackage{graphicx}

 \newcommand{\R}{\mathbb{R}}
 \newcommand{\SO}[1]{\mathrm{SO}_{#1}}
 \newcommand{\SE}[1]{\mathrm{SE}_{#1}}
 \newcommand{\trans}[1]{#1^{\mathsf T}}
 \newcommand{\B}{\mathrm{Bl}}
 \newcommand{\PJ}{\mathbb{P}}
 \newcommand{\x}{\underline{x}}
 \newcommand{\y}{\underline{y}}
 \newcommand{\w}{\underline{w}}
 
 \newtheorem{thm}{Theorem}

\begin{document}

\title[Compactification of the group of rigid motions]{Compactification of the group of rigid motions and applications
to robotics} 

\author{Nestor Djintelbe}
\address{N.~D.: Universit\'e Assane Seck, Ziguinchor, S\'en\'egal}
\thanks{N.~D. partially supported by EMS-Simons for Africa, PRMAIS network, Univ.~Rennes 1}

\author{Michel Coste}
\address{M.~C.: IRMAR, Universit\'e de Rennes I, Campus de Beaulieu, 35 042 Rennes, France}
\thanks{M.~C. partially supported by ANR KAPAMAT}


\begin{abstract}
We introduce a compactification of the group of rigid motions in 3-space derived from the Study model for this group. We use this compactification in robot kinematics, by considering the boundary of the configuration space of a robot. We study in particular the degeneration of the direct kinematic problem for some parallel robots with three degrees of freedom.
\end{abstract}

\maketitle

\section{Introduction}

The Lie group of rigid motions in 3-space is a fundamental tool in kinematics. This is a noncompact group, and therefore it is natural to look for compactifications. Such a compactification, known for a long time, is given by the Study model, as a quadric in a projective space of dimension 7. The Study model has proved very useful for the algebraic modelization of problems in robotics and their resolution in computer algebra systems. We give a presentation of the Study model in Section 2. In Section 3, we recall a few examples of the use of the Study model for algebraic modelization in robotics; this section also serves to introduce the parallel robots that we shall consider. This material is not new and we give references to the literature. 

In order to give a convenient kinematic interpretation to the boundary points, it is necessary to modify the compactification by blowing-up the Study quadric along its exceptional 3-plane, the complement of the image of the group of rigid motions. This is explained in Section 4, where we also show that this new compactification is isomorphic to the compactification obtained by taking the projective closure of the subgroup of translations. We describe geometrically the boundary of configuration space for the robots presented in the preceding section. The purpose of the introduction of boundary points is to obtain information on the kinematics of the robots as the lengths of their limbs are large enough. 

A prominent problem in the kinematics of parallel robots is the direct kinematic problem (DKP): what are the possible configurations of the robot for given actuated joint variables (for given lengths of limbs, in the cases under consideration). We propose in Section 5 to study a ``degenerate'' direct kinematic problem on the boundary of the configuration space of the robot. This degenerate DKP is much simpler than the full DKP of the robot, and it gives a rather faithful picture of the behaviour of the robot for sufficiently large limb lengths.

\section{Study model}
We present in this section the parametrization of the group of rigid motions in 3-space which was invented by E.~Study in 1891  \cite{Study91}. This parametrization is well adapted to robot kinematics (see \cite{HS10}) and has been used in many papers. Our presentation is a little different from the usual one.

\subsection{Rigid motions and rotations over dual numbers}
We denote by $\R[\epsilon]$ the algebra of dual numbers $a+\epsilon b$, where $a,b$ are real and $\epsilon^2=0$. An element of $\SO3(\R[\epsilon])$ is of the form $\mathbf R+\epsilon \mathbf M$, where $\mathbf R$ and $\mathbf M$ are real $3\times 3$ matrices such that 
$$\mathbf I_3=(\mathbf R+\epsilon\mathbf M)\trans{(\mathbf R+\epsilon \mathbf M)}=\mathbf R\,\trans{\mathbf R} + \epsilon(\mathbf R\,\trans{\mathbf M}+\mathbf M\,\trans{\mathbf R})\;.$$
 So $\mathbf R\in \SO3(\R)$ and $\mathbf M\,\trans{\mathbf R}$ is skew-symmetric. Every skew-symmetric $3\times 3$ matrix is the matrix of the cross-product $\mathbf x\mapsto \mathbf t\times \mathbf x$ for some vector $\mathbf t$; we denote this matrix by
$$ \Omega_{\mathbf t}=\begin{bmatrix} 0&-c&b\\ c&0&-a\\ -b&a&0\end{bmatrix}\quad\text{where } \mathbf t=\begin{bmatrix} a\\ b\\ c\end{bmatrix}\;.$$
Hence, every element of $\SO3(\R[\epsilon])$ can be uniquely written in the form: 
$$(\mathbf{I}_3+\epsilon \Omega_\mathbf{t}) \mathbf{R} \quad\text{where } \mathbf{R}\in \SO3(\R) \text{ and }\mathbf{t}\in \R^3\;.$$ 

To each element $(\mathbf{I}_3+\epsilon \Omega_\mathbf{t}) \mathbf{R} \in \SO3(\R[\epsilon])$ we associate the rigid motion $\mathbf{x}\mapsto \mathbf{R} \mathbf{x}+\mathbf{t}$ in $\SE3(\R)$. We denote by $\Phi$ this bijection.

\begin{thm}
$\Phi : \SO3(\R[\epsilon])\to \SE3(\R)$ is a group isomorphism.
\end{thm} 

\noindent \textit{Proof.} Let $(\mathbf{I}_3+\epsilon \Omega_\mathbf{t}) \mathbf{R}$ and $(\mathbf{I}_3+\epsilon \Omega_\mathbf{v}) \mathbf{S}$ be two elements of $\SO3(\R[\epsilon])$.  
On the one hand we have
$$\mathbf{R}(\mathbf{S} \mathbf{x}+\mathbf{u})+\mathbf{t}=(\mathbf{R} \mathbf{S})\mathbf{x}+(\mathbf{R} \mathbf{u} +\mathbf{t})\;,$$
and on the other hand 
$$\begin{aligned}(\mathbf{I}_3+\epsilon \Omega_\mathbf{t}) \mathbf{R} (\mathbf{I}_3+\epsilon \Omega_\mathbf{u}) \mathbf{S} &= \mathbf{R}\mathbf{S}+\epsilon(\Omega_\mathbf{t} \mathbf{R}\mathbf{S} + \mathbf{R}\Omega_\mathbf{u} \mathbf{S})\\ &=\mathbf{R}\mathbf{S}+\epsilon(\Omega_\mathbf{t} \mathbf{R}\mathbf{S} +\Omega_{\mathbf{R} \mathbf{u}}\mathbf{R} \mathbf{S})=(\mathbf{I}_3+\epsilon\Omega_{\mathbf{R} \mathbf{u}+\mathbf{t}})\mathbf{R} \mathbf{S}\;.
\end{aligned}$$
This proves that the bijection $\Phi$ is indeed a group isomorphism. \hfill $\square$
\par\medskip

\noindent \textbf{Remark.}
A point with values in $\R[\epsilon]$ is the same as a real point plus a tangent vector at this point. Moreover, the Lie algebra of $\SO3(\R)$ may be identified with $\R^3$ equipped with the cross-product and the adjoint action, with this identification, is the action by rotations. This is what is behind Theorem~1.

\subsection{Dual quaternions}
It is well-known that the group of unit quaternions (quaternions with norm 1) is a double covering of $\SO3$. 
The rotation matrix image of the quaternion $x_0+x_1\mathbf{i}+x_2\mathbf{j}+x_3\mathbf{k}$ is
\begin{equation}\label{eq:matrot}
\begin{aligned}
&\mathbf{\mathbf{R}}(x_0,x_1,x_2,x_3) =\\
&\quad \begin{bmatrix} x_0^2 + x_1^2 - x_2^2 - x_3^2 & 2(x_1x_2 - x_0x_3) & 2(x_1x_3 + x_0x_2)\\ 2(x_1x_2 + x_0x_3) & x_0^2 - x_1^2 + x_2^2 - x_3^2 & 2(x_2x_3 - x_0x_1)\\ 2(x_1x_3 - x_0x_2) & 2(x_2x_3 + x_0x_1) & x_0^2 - x_1^2 - x_2^2 + x_3^2 \end{bmatrix}
\end{aligned}
\end{equation}
This is also the case working over $\R[\epsilon]$ instead of $\R$. The unit elements in the algebra $\mathbb{H}[\epsilon]$ of dual quaternions are those $\mathbf{q}+\epsilon \mathbf{r}$ with $\mathbf{q},\mathbf{r}\in \mathbb{H}$ such that 
$$1=(\mathbf{q}+\epsilon \mathbf{r})(\overline{\mathbf{q}}+\epsilon \overline{\mathbf{r}})= \mathbf{q} \overline{\mathbf{q}} + \epsilon (\mathbf{q}\overline{\mathbf{r}}+ \mathbf{r}\overline{\mathbf{q}})\;.$$
With $\mathbf{q}=x_0+x_1\mathbf{i}+x_2\mathbf{j}+x_3\mathbf{k}$ and $\mathbf{r}=y_0+y_1\mathbf{i}+y_2\mathbf{j}+y_3\mathbf{k}$, this amounts to
$$ x_0^2+x_1^2+x_2^2+x_3^2 =1 \quad \text{and}\quad x_0y_0 + x_1y_1 + x_2y_2 + x_3y_3 = 0\;.$$

The element of $\mathrm{SO}_3(\R[\epsilon])$ image of the unit dual quaternion $\mathbf{q}+\epsilon \mathbf{r}$ in the double covering is $(\mathbf{I}_3+\epsilon\Omega_\mathbf{t})\mathbf{R}$ where $\mathbf{R}$ is given by formula (\ref{eq:matrot}) and
\begin{equation}\label{eq:vecttrans}
\mathbf{t}(x_0,x_1,x_2,x_3,y_0,y_1,y_2,y_3)=2\,\begin{bmatrix} 
x_0y_1 - x_1y_0 + x_2y_3 - x_3y_2\\
x_0y_2 - x_1y_3 - x_2y_0 + x_3y_1\\
x_0y_3 + x_1y_2 - x_2y_1 - x_3y_0
 \end{bmatrix}
\end{equation}
This formula (\ref{eq:vecttrans}) is obtained by plugging in the dual quaternion $(\mathbf{q}+\epsilon\mathbf{r})\overline{\mathbf{q}}=1+\epsilon\mathbf{r}\overline{\mathbf{q}}$ in formula (\ref{eq:matrot}) and identifying the result with $\mathbf{I}_3+\epsilon\Omega_\mathbf{t}$; note that $\mathbf{r}\overline{\mathbf{q}}$ is a pure quaternion.

\subsection{Study quadric}
Formulas (\ref{eq:matrot}) and (\ref{eq:vecttrans}) are homogeneous of degree 2 in variables $\x,\y$ (we use $\x$ for $x_0,\ldots,x_3$ and $\y$ for $y_0,\ldots,y_3$). So instead of describing rigid motions by pairs of opposite unit dual quaternions, we can use points in the 7-dimensional real projective space $\PJ^7(\R)$ with homogeneous coordinates $\x,\y$ satisfying $\sum_{i=0}^3 x_iy_i=0$ and $\sum_{i=0}^3 x_i^2\neq 0$. In this way we arrive to the Study model for rigid motions: 

\begin{thm}[Study]
Let $S\subset \PJ^7(\R)$ be the 6-dimensional \emph{Study quadric} with equation 
$$\sum_{i=0}^3 x_iy_i=0\;,$$ 
and let $E\subset S$ be the 3-plane contained in $S$ with equations  $x_0=x_1=x_2=x_3=0$.

There is a one-to-one correspondance which, to each point with homogeneous coordinates $\x,\y$ belonging $S\setminus E$, associates the rigid motion
$$ \mathbf{u}  \longmapsto  \dfrac{1}{\Delta(\x)}\left(\mathbf{R}(\x)\mathbf{u} +\mathbf{t}(\x,\y)\right)\;,$$
where $\Delta(\x)=\sum_{i=0}^3x_i^2$, $\mathbf{R}$ is given by formula (\ref{eq:matrot}) and $\mathbf{t}$  by formula (\ref{eq:vecttrans}). \end{thm}

The Study quadric contains 3-planes. There are two 6-dimensional families of these 3-planes, and these 3-planes have kinematic significance. For instance the 3-plane $y_0=y_1=y_2=y_3=0$ (corresponding to rotations with center the origin) is in one family which also contains the ``exceptional'' 3-plane  $E$; the 3-plane $ y_0=x_1=x_2=x_3=0$ (corresponding to translations) is in the other family. We shall see others of these 3-planes in the next section. For more details we refer to \cite{Sel05} p.~246.

\section{Operation modes of parallel robots}

\subsection{Configuration space and modes of operation}

We report in this section a few results concerning the kinematics of mobile platforms with different architectures which are estabished in the papers \cite{SWH14, WHP09, WH01}
We are considering robots consisting of a mobile platform linked to a fixed base by three limbs with variable lengths (their lengths are controlled by actuated prismatic (P) joints). Both base and platform are equilateral triangles, and the limbs are attached to the base and the platform via joints centered at the vertices $A_i$ and $B_i$ of these triangles. The mobile platform has restricted degrees of freedom due to the specification of the joints.

For the algebraic modelization we work in the fixed frame attached to the base:

\begin{center}
\includegraphics[scale=.4]{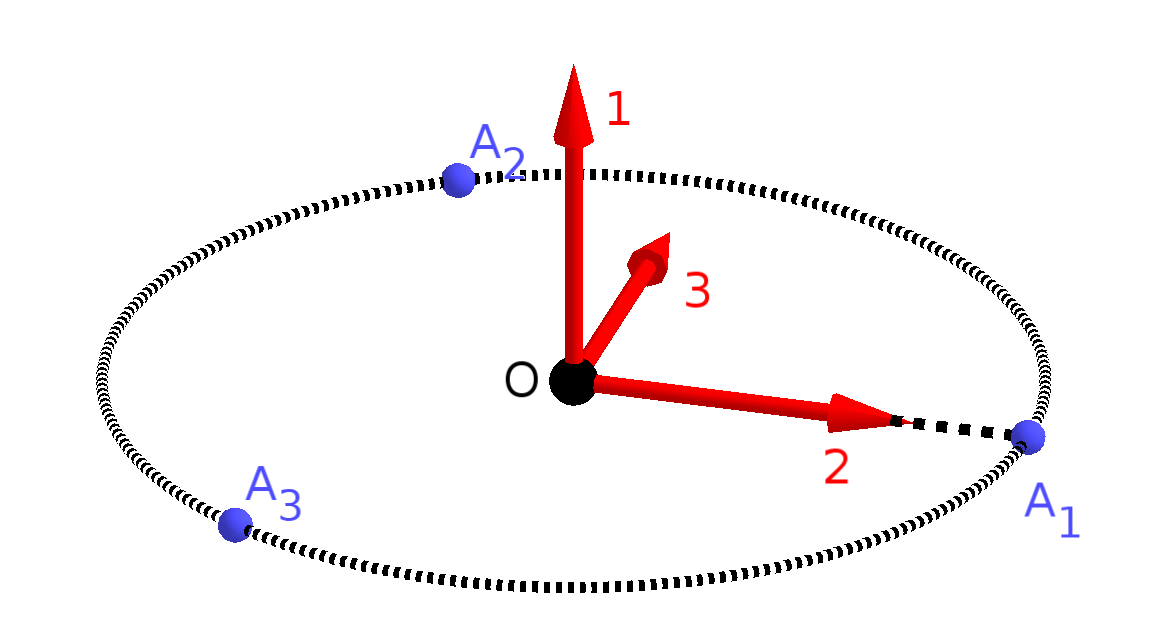}
\end{center}

In this frame the vertices of the base have coordinates
$$A_1: \begin{bmatrix}0\\ k_1\\0\end{bmatrix}\quad A_2: \begin{bmatrix}0\\-k_1/2\\ \sqrt3\, k_1/2\end{bmatrix}\quad A_3: \begin{bmatrix}0\\ -k_1/2\\ -\sqrt3\, k_1/2\end{bmatrix}\;,$$
where $k_1$ is the radius of the circle circumscribed to the base. We have a similar mobile frame attached to the platform (with $k_2$ the radius of its circumscribed circle) and we compute the coordinates of the points $B_i$ in the fixed frame using the Study parametrization with homogeneous coordinates $(\x,\y)$ for the change of frame. If we denote by $b_i$ the vector of coordinates of $B_i$ in the frame attached to the platform, we have
$$B_i= \dfrac1{\Delta(\x)}(\mathbf{R}(\x) b_i + \mathbf{t}(\x,\y))\;.$$

We form the homogeneous ideal generated by the constraint equations of the platform and  the equation of the Study quadric $S$. We saturate this ideal with respect to $\Delta(x)=\sum_{i=0}^3 x_i^2$ in order to remove spurious components contained in the exceptional 3-plane $E$. The ideal obtained is the ideal of the \emph{configuration space} of the mobile platform, an algebraic subset of the Study quadric.This configuration space is not irreducible in the examples we shall consider. Its irreducible components are the \emph{modes of operation} of the platform. We obtain these different modes of operation by computing the primary decomposition of the ideal of the configuration space.
In the examples, we shall perform the computations for $k_1=1$ and $k_2=3/2$.

\subsection{3-RPS \cite{SWH14}}\label{bord3RPS}

In the 3-RPS architecture, each limb is attached to the base via a passive revolute (R) joint whose axis of rotation is tangent to the circle circumbscribed to the platform, and attached to the platform via a passive spherical (S) joint which allows all rotations around its center. The limb $(A_iB_i)$ is perpendicular to the axis of rotation at $A_i$.
 
 \begin{center}
   \includegraphics[scale=.27]{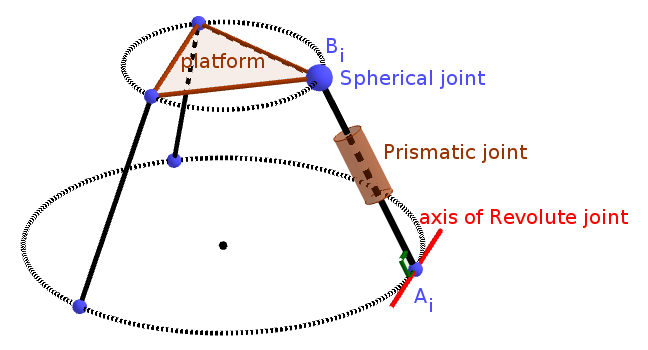}
 \end{center}
 
 The constraint equation for each limb expresses that $(A_iB_i)$ is orthogonal to the tangent at $A_i$ to the circumscribe circle. The configuration space is 3-dimensional and there are two modes of operation whose ideals are:
 \begin{itemize}
  \item $\mathfrak I_1 = \langle x_0 ,\; x_1 y_1 + x_2 y_2 + x_3 y_3 ,\; 3 x_2 x_3 + 2 x_2 y_0 - 2 x_3 y_1 + 2 x_1 y_3 ,\; 3 x_2^2 - 3 x_3^2 + 4 x_3 y_0 + 4 x_2 y_1 - 4 x_1 y_2 ,\; 9 x_3^2 y_1 - 6 x_3 y_0 y_1 - 8 y_0^2 y_1 - 6 x_2 y_1^2 - 8 y_1^3 - 8 y_1 y_2^2 - 18 x_3 y_2 y_3 - 6 x_2 y_3^2 - 8 y_1 y_3^2 ,\; 9 x_3^3 - 6 x_3^2 y_0 - 8 x_3 y_0^2 + 4 x_2 y_0 y_1 - 12 x_3 y_1^2 + 12 x_1 x_3 y_2 + 8 x_1 y_0 y_2 + 6 x_1 x_2 y_3 - 12 x_2 y_2 y_3 - 12 x_3 y_3^2 \rangle$
 \item $\mathfrak I_2 = \langle x_1 ,\; x_0 y_0 + x_2 y_2 + x_3 y_3 ,\; 3 x_2 x_3 + 2 x_2 y_0 - 2 x_3 y_1 - 2 x_0 y_2 ,\; 3 x_2^2 - 3 x_3^2 + 4 x_3 y_0 + 4 x_2 y_1 - 4 x_0 y_3 ,\; 9 x_3^2 y_0 - 6 x_3 y_0^2 - 8 y_0^3 - 6 x_2 y_0 y_1 - 8 y_0 y_1^2 + 6 x_3 y_2^2 - 8 y_0 y_2^2 + 6 x_2 y_2 y_3 - 12 x_3 y_3^2 - 8 y_0 y_3^2 ,\; 27 x_3^3 - 36 x_3 y_0^2 - 16 y_0^3 - 36 x_3 y_1^2 - 16 y_0 y_1^2 - 18 x_0 x_2 y_2 - 36 x_0 y_1 y_2 + 12 x_3 y_2^2 - 16 y_0 y_2^2 + 36 x_0 x_3 y_3 - 12 x_2 y_2 y_3 - 48 x_3 y_3^2 - 16 y_0 y_3^2 \rangle$
 \end{itemize}
 The lists of generators of these ideals is not very informative. Nevertheless we can see that in the first mode, the rotation part is a half-turn ($x_0=0$) while the rotation part in the second mode has an horizontal axis ($x_1=0$). One can also check (although this is not visible on the computed generators) that the two modes are interchanged by the involution consisting in multiplying on the right with quaternion $\mathbf i$, that is composing on the right with the half-turn with vertical axis through the origin.

\subsection{SNU 3-UPU \cite{WHP09}}

In the 3-UPU architectures, each limb is attached to the base and to the platform via passive Cardan joints, also called universal joints (U). Each U joint has two axes of rotation intersecting orhogonally in the center of the joint. We shall study two different architectures. In the SNU 3-UPU architecture, the rotation axis of the U-joint which is rigidly fixed on the base (resp. platform) is pointing towards the center of its circumscribed circle. The two rotation axes rigidly fixed on the limb are parallel, and both orthogonal to the limb. In the following picture, the rotation axes are numbered following the kinematice chain from base to platform.

 \begin{center}
   \includegraphics[scale=.52]{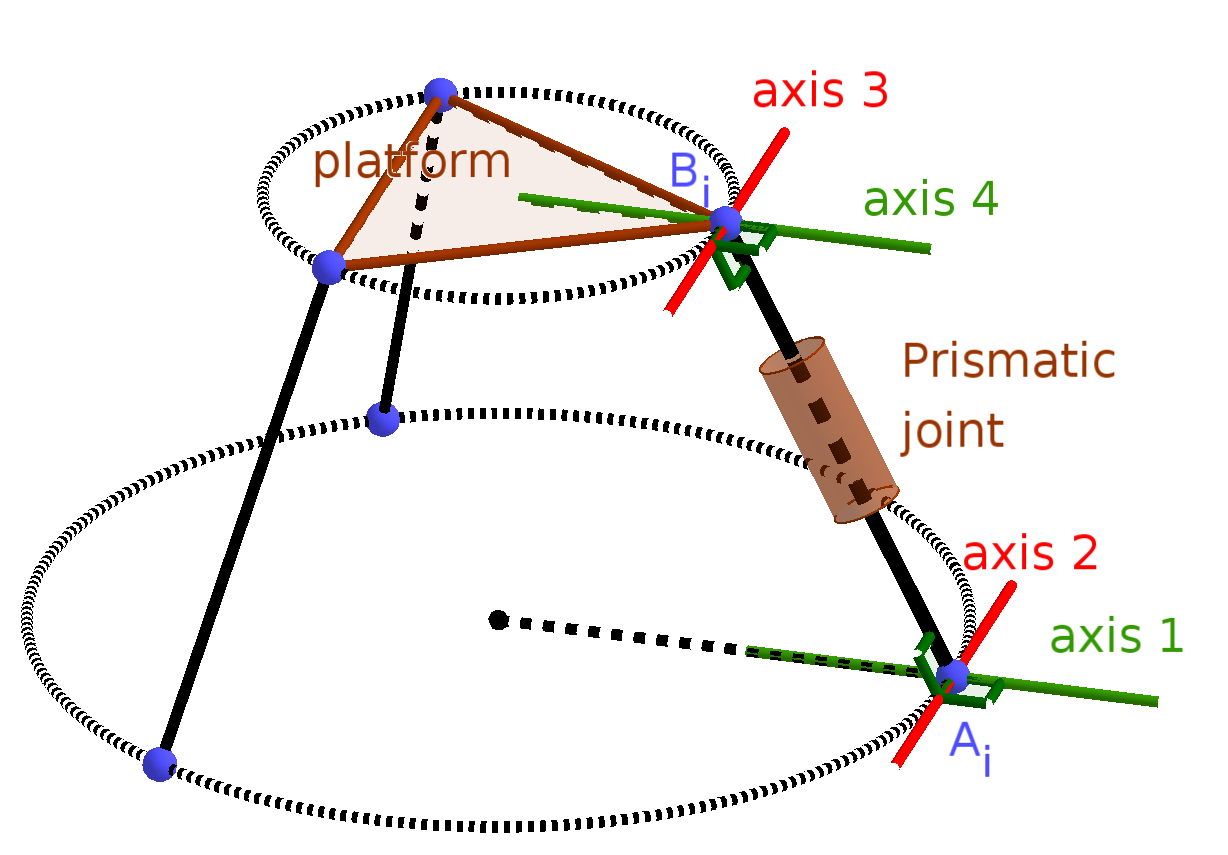}
 \end{center}

The constraint equation for each limb expresses that axis 1, axis 4 and $(A_iB_i)$ are coplanar. There are seven  modes of operation with 3 d.o.f., all 3-planes in the Study quadric, whose ideals are :
\begin{itemize}
\item $\mathfrak K_0=\langle y_0,y_1,y_2,y_3\rangle$ ; all rotations around the origin.
\item $\mathfrak K_1=\langle y_0,x_1,x_2,x_3\rangle$ : all translations.
\item $\mathfrak K_2=\langle x_0,y_1,x_2,x_3\rangle$ : rigid motions consisting in the half-turn with vertical axis through origin, followed by a translation.
\item $\mathfrak K_3=\langle y_0,y_1,x_2,x_3\rangle$ : all rigid motions in the base plane.
\item $\mathfrak K_4=\langle x_0,x_1,y_2,y_3\rangle$ : rigid motions consisting in a horizontal flip followed by motion in the base plane.
\item $\mathfrak K_5=\langle x_0,y_1,y_2,y_3\rangle$ : rigid motions consisting in a half-turn with axis through origin, followed by a translation in the direction of the axis of half-turn.
\item $\mathfrak K_6=\langle y_0,x_1,y_2,y_3\rangle$ : rigid motions consisting of the half-turn with vertical axis through origin, followed by a rigid motion of mode $\mathfrak K_5$. The rotation part of a rigid motion in this mode is a rotation with horizontal axis.
\end{itemize}
plus a non-real component $\mathfrak K_7$ whose real points correspond to translations along the vertical axis through origin, possibly composed with the half-turn around this axis. These are contained in other modes of operation.

The configuration space is stable by the involution consisting in composing on the right with the half-turn with vertical axis through the origin. Modes $\mathfrak K_0, \mathfrak{K}_3, \mathfrak{K}_4$ and $\mathfrak K_7$ are stable by this involution which exchanges modes $\mathfrak{K}_1, \mathfrak{K}_2$ and also modes $\mathfrak{K}_5, \mathfrak{K}_6$.

\subsection{Tsai 3-UPU \cite{WH01}}

The only thing which makes the Tsai 3-UPU architecture different from the SNU one is the fact that rotation axes 1 and 4 on each limb are tangent to the circumscribed circles of base and platform respectively. 

 \begin{center}
  
\includegraphics[scale=.52]{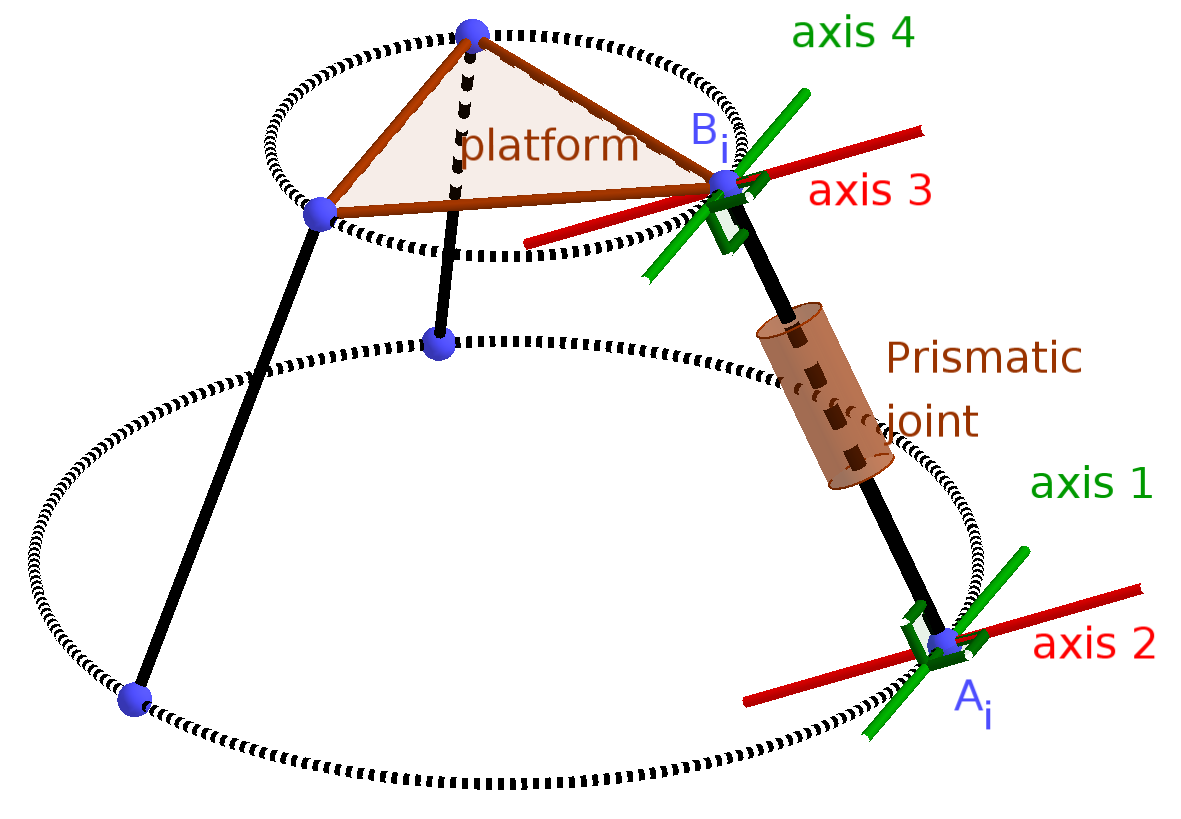}
 \end{center}

This small change changes drastically the kinematic behaviour of the platform. The configuration space now contains four modes of operation corresponding to 3-planes in the Study quadric, which are the modes $\mathfrak{K}_1,\mathfrak{K}_2,\mathfrak{K}_3,\mathfrak{K}_4$ of the SNU, plus another ``mysterious'' mode which is a real subvariety of dimension 3 of the Study quadric whose computed ideal has a long list of generators:
$$\begin{aligned}
\mathfrak K_8&=\langle 2x_1y_1 + x_2y_2 + x_3y_3, 2x_0y_0 + x_2y_2 + x_3y_3, \ldots \\
&\qquad \ldots, 15x_2^2x_3y_0 - 5x_3^3y_0 - x_2^3y_1 + 3x_2x_3^2y_1 - x_2^2y_0y_1 - x_3^2y_0y_1\rangle
 \end{aligned}$$
The kinematic analysis of this last mode of operation is not easy. We shall see later how to recover partial information on this mode, via degeneration.

\section{Compactification and boundary}

We introduce in this section a compactification of the group of rigid motions, and also compactifications of configuration spaces and modes of operation. The purpose of these compactifications is to use degeneration to the boundary in order to obtain some partial information on the kinematic behaviour of the robots. The usefulness of points in the boundary for kinematic problems is exemplified in \cite{COSTE201582} and also, for a different kind of boundary points, in \cite{HEGEDUS2013407}.

\subsection{Blowing-up of $S$ along $E$}

The 3-plane $E$ with equations $x_0=x_1=x_2=x_3=0$ may be seen as the boundary of $\SE3(\R)$ in the Study quadric $S$. This boundary is composed of limits of rigid motions as the norm of the translation vector goes to infinity. In this boundary, all information about the rotation part of the rigid motion is lost. 

To avoid this, we blow up $S$ along $E$ in order to replace $E$ with a variety of dimension 5 while $E$ is of dimension 3 only. The information on rotation part will not be lost in this new boundary.  Blowing-up a variety $V$ gives a new variety that is birational to $V$; some information on blowing-up is given in \cite{BCR98}, p.~78.

In order to describe more precisely the blowing-up, we introduce new homogeneous variables $\w=w_0,w_1,w_2,w_3$. The blowing-up of $S$ along $E$ is $\B_E(S)$ with the projection $
\pi : \B_E(S) \to S$
where
\begin{multline}\label{eq: blowup}
\B_E(S)=\Big\{ \left(\left[\x,\y\right],\left[\w\right]\right) \in \PJ^7(\R)\times \PJ^3(\R) \mid  \sum_{i=0}^3 x_i y_i=0,\quad \sum_{i=0}^3 y_i w_i=0 \\  x_i w_j-x_j w_i=0\quad\text{for } 0\leqslant i<j\leqslant 3
 \Big\}
\end{multline}
and $\pi$ is the restriction of the projection on the first factor.

$\B_E(S)$ is a compact variety of dimension 6. The restriction of $\pi$ to $\B_E(S) \setminus \pi^{-1}(E)$ is an isomorphism onto $S \setminus E$, so that $\SE3(\R)$ may be identified with $\B_E(S) \setminus \pi^{-1}(E)$ and $\pi^{-1}(E)$  appears as the boundary of $\SE3(\R)$ in $\B_E(S)$. We denote by 
\[
\alpha:\mathrm{SE3(\R)}\cong S\setminus E\longrightarrow   \B_E(S) \setminus \pi^{-1}(E)\hookrightarrow \B_E(S)
\]
the compactification of $\SE3(\R)$ thus obtained. 

The inverse image $\pi^{-1}(E)$ of $E$ in the blowing-up, i.e. \emph{the exceptional divisor}, is the hypersurface in the product of projective spaces $\PJ^3(\R)\times \PJ^3(\R)$ (with homogeneous coordinates $\y$ for the first factor and $\w$ for the second) whose equation is
\[
w_0y_0+w_1y_1+w_2y_2+w_3y_3=0\;.
\]
It looks like the equation for the Study quadric, but here it is an equation of bidegree $(1,1)$ which describes a variety of dimension $5$. The homogeneous coordinates $ \w$ and $\y$ give respectively information about the rotation (in the usual way that a quaternion determines a rotation) and the direction of the translation vector by the formula (\ref{eq:vecttrans}) by replacing the $ \x $ with the corresponding $ \w $. More precisely this direction is given by the vector 
\[ \mathbf u (\w,\y) = \begin{bmatrix} 
w_0y_1 - w_1y_0 + w_2y_3 - w_3y_2\\
w_0y_2 - w_1y_3 - w_2y_0 + w_3y_1\\
w_0y_3 + w_1y_2 - w_2y_1 - w_3y_0 \end{bmatrix}\;.\] 
If we normalize the homogeneous coordinates $\w$ by $\sum_{i=0}^3w_i^2=1$ and $\y$ by $\sum_{i=0}^3y_i^2=1$, then the vector $\mathbf{u}(\w,\y)$ has norm $1$.

\subsection{Comparison of compactifications}

The group $\SE3(\R)$ is the semi-direct product  $\SO3(\R)  \ltimes \R^3$. The rotation group $\SO3(\R)$ is compact, and it is isomorphis as real variety to $\PJ^3(\R)$ via the regular mapping $[\x] \mapsto \frac{1}{\Delta(\x)}  \mathbf{R}(\x)$ where $\mathbf{R}$ is as in (\ref{eq:matrot}). We can compactify $\R^3$ to $\PJ^3(\R)$.  Thus we obtain the following compactification of $\SE3(\R)$:
$$
\beta:\SE3(\R) \cong \SO3(\R) \ltimes \R^3 \hookrightarrow \PJ^3(\R) \times \PJ^3(\R)\;.
$$
We are going to see that this compactification coincides with the one obtained by blowing-up the Study quadric.
Define the morphism $\sigma$
$$
\begin{aligned}
\sigma: \PJ^3(\R) \times \PJ^3(\R) & \longrightarrow  \PJ^7(\R) \times \PJ^3(\R)\\ 
 ([\w],[r,s,t,u]) & \longmapsto ([\x,\y],[\w])
\end{aligned}
$$
by
\begin{equation}\label{eq: formulessigma}
\left \{ \begin{array}{lcl}
x_0=w_0r &\qquad& 2y_0=  -w_1s - w_2t - w_3u\\
x_1=w_1r && 2y_1=\phantom{-} w_0s+w_3t-w_2u\\
x_2=w_2r && 2y_2= -w_3s+w_0t+w_1u\\
x_3=w_3r && 2y_3= \phantom{-}w_2s-w_1t+w_0u
\end{array}\right.\;.
\end{equation}
\begin{thm}\label{th: isomcompact}
The two compactifications $\alpha$ and $\beta$ of $\SE3(\R)$ are isomorphic. More precisely, the morphism $\sigma$ induces a biregular isomorphism of real varieties (in the sense of \cite{BCR98}) $\sigma: \PJ^3(\R) \times \PJ^3(\R)\to\B_E(S)$ such that $\alpha=\sigma\circ \beta$.
\end{thm} 

\noindent \textit{Proof}.
The mapping $\sigma$ is a well-defined regular mapping. Indeed, the formulas (\ref{eq: formulessigma}) are homogeneous of bidegree $(1,1)$ and, since 
\begin{equation}
\begin{aligned}\label{eq: formulesinv}
(w_0^2+w_1^2+w_2^2+w_3^2) r &= w_0x_0 + w_1x_1 + w_2x_2 + w_3x_3 \\
(w_0^2+w_1^2+w_2^2+w_3^2) s &=2\left( - w_1y_0+ w_0y_1 - w_3y_2 + w_2y_3\right)\\
(w_0^2+w_1^2+w_2^2+w_3^2) t &=2\left( - w_2y_0  + w_3y_1 +w_0y_2- w_1y_3\right)\\
(w_0^2+w_1^2+w_2^2+w_3^2) u &=2\left( - w_3y_0 - w_2y_1 + w_1y_2+w_0y_3\right)\;,
\end{aligned}
\end{equation}
if all $x_i$ and all $y_i$ are zero, then all $w_i$ are zero or $r=s=t=u=0$, which is impossible for homogeneous coordinates. It is also clear from formulas (\ref{eq: formulessigma}) and the description (\ref{eq: blowup}) of $\B_E(S)$ that the image of $\sigma$ is contained in $\B_E(S)$. In order to show that $\sigma$ is an isomorphism of real algebraic varieties from $\PJ^3(\R)\times \PJ^3(\R)$ onto $\B_E(S)$, we show that the regular mapping
\[
\begin{aligned}
\tau:  \PJ^7(\R) \times \PJ^3(\R)\supset \B_E(S)& \longrightarrow \PJ^3(\R) \times \PJ^3(\R),\\ 
([\x,\y],[\w])& \longmapsto ([\w],[r,s,t,u]),
\end{aligned}
\]
defined by the homogeneous formulas of bidegree $(1,1)$
\begin{equation}
\begin{aligned}\label{eq: formulesipi}
r &= w_0x_0 + w_1x_1 + w_2x_2 + w_3x_3 \\
s &=2\left( - w_1y_0+ w_0y_1 - w_3y_2 + w_2y_3\right)\\
t &=2\left( - w_2y_0  + w_3y_1 +w_0y_2- w_1y_3\right)\\
u &=2\left( - w_3y_0 - w_2y_1 + w_1y_2+w_0y_3\right)\;,
\end{aligned}
\end{equation}
is the inverse of $\sigma$. 
The morphism $\tau$ is well-defined on $\B_E(S)$. Indeed the $r,s,t,u$ given by formulas (\ref{eq: formulesipi}) cannot vanish simultaneaously on $\B_E(S)$.
From formulas (\ref{eq: formulesipi}) and $w_ix_j=w_jx_i$, we obtain:
\begin{equation}\label{eq: formulesinvbis}
w_i r = (w_0^2+w_1^2+w_2^2+w_3^2) x_i\qquad \text{for } i=0,\ldots, 3
\end{equation}
Hence, if $r=0$ then $x_i=0$ for $i=0,\ldots, 3$.

Also, from formulas (\ref{eq: formulesipi}) and $\sum_{i=0}^3 w_iy_i=0$ we obtain:
\begin{equation}\label{eq: formulesinvter}
\begin{aligned}
  -w_1s - w_2t - w_3u &= 2(w_0^2+w_1^2+w_2^2+w_3^2)y_0\\
\phantom{-} w_0s+w_3t-w_2u &= 2(w_0^2+w_1^2+w_2^2+w_3^2)y_1\\
 -w_3s+w_0t+w_1u &= 2(w_0^2+w_1^2+w_2^2+w_3^2)y_2\\
 \phantom{-}w_2s-w_1t+w_0u &= 2(w_0^2+w_1^2+w_2^2+w_3^2)y_3
\end{aligned}
\end{equation}
Hence, if $s=t=u=0$, then $x_i=0$ for $i=0,\ldots, 3$. As $x_i$ and $y_i$ cannot all be $0$, this show that $\tau$ is well-defined. 

Furthermore, from equations (\ref{eq: formulesinv}) we see that $\tau\circ \sigma$ is the identity of $\PJ^3(\R) \times \PJ^3(\R)$, while equations (\ref{eq: formulesinvbis}) and (\ref{eq: formulesinvter}) show that $\sigma\circ \tau$ is the identity of $\B_E(S)$.

It remains to show that $\alpha=\sigma\circ \beta$, which is equivalent to $\beta=\tau\circ \alpha$ ; the latter equality is clear by comparing formulas (\ref{eq:vecttrans}) for the translation vector in terms of Study parameters and formulas (\ref{eq: formulesipi}) for $\tau$, taking into account the fact that the $\x$ and the $\w$ are proportional on $\alpha(\SE3(\R))$. This concludes the proof of Theorem \ref{th: isomcompact}.\hfill$\square$
\par\medskip

In the following we shall use the compactification $\alpha$, since it is the one which is well adapted to Study parameters, and Study parameters are very convenient for computations in kinematics.

\subsection{Boundary of configuration space}\label{sec: bdryconf}

Let $C$ be a configuration space or a mode of operation of a mechanism, identified with an algebraic subset of the Study quadric $S$. We can take its strict transform $\widetilde{C}$ in the blowing up $\pi : \B_E(S)\to S$: this is the inverse image $\pi^{-1}(C)$ with all components contained in $\pi^{-1}(E)$ removed. Equations for $\widetilde{C}$ can be computed by adding to the ideal of $C$ the equations of $\B_E(S)$ and saturating with respect to $\sum_{i=0}^3 x_i^2$.
The boundary of $C$ is then $\tilde{C} \cap \pi^{-1}(E)$. Setting $\x=0$ in the equations for $\widetilde{C}$, we obtain bihomogeneous equations in $\w,\y$ for the boundary of $C$ in $\PJ^3(\R)\times \PJ^3(\R)$. We give now two examples of this computation of boundary

\subsubsection{Boundaries of the modes of operation of the 3-RPS}

The boundary of the mode of operation $\mathfrak{I}_1$ is the algebraic subset $\mathfrak{J}_1$ of $\PJ^3(\R)\times \PJ^3(\R)$ with equations 
\[
w_0=y_1=w_2y_2+w_3y_3=w_1y_3+w_2y_0=w_1y_2-w_3y_0=0\;.
\]
This algebraic subset is a projective 2-plane, we can choose $w_1,w_2,w_3$ as homogeneous coordinates for this 2-plane. We can then take $y_0= w_1,y_2=w_3, y_3=-w_2$. The points in this boundary may be seen as rigid motions whose rotation part is a half-turn and translation part infinite translation in the vertical direction. \par\medskip
  
The boundary of the mode of operation $\mathfrak{I}_2$ is the algebraic subset $\mathfrak{J}_2$ of $\PJ^3(\R)\times \PJ^3(\R)$ with equations 
\[
w_1=y_0=w_2y_2+w_3y_3=w_0y_3-w_2y_1=w_0y_2+w_3y_1=0\;.
\]
This is also a projective 2-plane, we can choose here $w_0,w_2,w_3$ as homogeneous coordinates for this 2-plane. We can then take $y_1= w_0,y_2=-w_3, y_3=w_2$. The points of $\mathfrak{J}_2$ may be seen as rigid motions whose rotation part is a rotation with horizontal axis and translation part infinite translation in the vertical direction. \par\medskip

\subsubsection{Boundary of the mysterious mode of operation of the Tsai 3-UPU}

The boundary of the mode of operation $\mathfrak{K}_8$ is not an irreducible algebraic subset of $\PJ^3(\R)\times \PJ^3(\R)$, but it decomposes into three components.

\begin{itemize}
\item The algebraic subset $\mathfrak L_{5}$ with equations $w_0=y_1=y_2=y_3=0$. This projective 2-plane is actually also the boundary of the mode of operation $\mathfrak{K}_5$ of the SNU 3-UPU. Its elements may be seen as rigid motions whose rotation part is a half-turn and translation part an infinite translation parallel to the axis of the half-turn
\item The algebraic subset $\mathfrak L_{6}$ with equations $y_0=w_1=y_2=y_3=0$, which is another 2-plane and the boundary of the mode of operation $\mathfrak{K}_6$ of the SNU 3-UPU. Its elements may be seen as half-turn with vertical axis followed by a ``rigid motion'' of $\mathfrak L_{6}$
\item A non-real component $\mathfrak{L}_7$ with equations 
$$\begin{aligned}
w_2^2 + w_3^2&=y_3w_2 + y_2w_3=y_2w_2 - y_3w_3=\\
&=y_1w_1 + y_3w_3=y_0w_0 + y_3w_3=y_2^2 + y_3^2=0,
\end{aligned}$$
which is the boundary of the component $\mathfrak{K}_7$ for the SNU 3-UPU and whose real points are only two singular points,  one in $\mathfrak L_{5}$ and the other in $\mathfrak L_{6}$
\end{itemize}

We see in this example how the degeneration to the boundary may give some information on the kinematic behavior of a mode of operation which is hard to analyze: we know how rigid motions in this mode of operation look like, when the lengths of the limbs become larger and larger. 

\section{Degeneration of the direct kinematic problem}

The direct kinematic problem (DKP) for the platforms we consider is the following: given the lengths $r_1,r_2,r_3$ of the limbs, what are the possible configurations for the platform? The lengths $r_i$ of the limbs are controlled by the actuated prismatic joints: they are the actuated joint variables. They can be easily computed from the homogeneous variables $\x,\y$ of points of the Study quadric. The mapping which associates to a point of the configuration space of the platform the triple $(r_1,r_2,r_3)\in \R^3$ is called the inverse kinematic mapping (IKM). The number of solutions to the DKP can change at the critical values of the IKM, the images of the singular points of the IKM. These singularities play a prominent role in the kinematics of parallel robots.

\subsection{Degeneration of the DKP}

We consider a degeneration of the DKP when the lengths of the limbs $(r_1,r_2,r_3)$ tend to infinity. We have already the boundary of the configuration space in $\PJ^3(\R)\times\PJ^3(\R)$. We need also to choose a boundary for the space of actuated joint variables. We cannot simply take $\PJ^3(\R)$ as compactification for the space of actuated joint variables since all points with homogeneous coordinates $(r_1,r_2,r_3,1)$ corresponding to configurations of the platform would tend to the point with homogeneous coordinates $(1,1,1,0)$ as the lengths of the limbs tend to infinity. Instead we make the change of variables $d_1=r_1-r_3$, $d_2=r_2-r_3$; we then let $r_3$ tend to infinity. More precisely, we embed $\R^2\times \R$ with variables $(d_1,d_2,r_3)$ into $\R^2\times \PJ^1(\R)$ and we take the plane $\R^2\times\{\infty\}$ as boundary for the space of actuated variables. 

We have now to describe the degenerate IKM from the boundary of configuration space (with bihomogeneous variables $\w,\y$) to the boundary of the space of actuated variables  (with variables $d_1,d_2$). As the limbs become parallel to the the translation vector as their lengths tend to infinity, the limit of the differences of lengths of limbs are, up to sign, differences of scalar products of the $\mathbf{R}(\w)b_i-A_i$ 
with the unit vector $\mathbf{u}(\w,\y)$ giving the direction of the infinite translation:
\begin{equation}\label{eq:degIKM}\begin{aligned}
d_1&= (\mathbf{R}(\w) b_1 - A_1)\cdot \mathbf{u}(\w,\y) - (\mathbf{R}(\w) b_3 - A_3)\cdot \mathbf{u}(\w,\y)\\
d_2&= (\mathbf{R}(\w) b_2 - A_2)\cdot \mathbf{u}(\w,\y) - (\mathbf{R}(\w) b_3 - A_3)\cdot \mathbf{u}(\w,\y)
\end{aligned}\end{equation}
We have to take into account in the discussion of the degenerate DKP that the signs of both $d_1$ and $d_2$ in (\ref{eq:degIKM}) depend on the orientation pf $\mathbf{u}$.

The degenerate DKP for a planar parallel robot with three degrees of freedom (3-RPR) has been considered in  \cite{Cos12}.

\subsection{Degenerate DKP for a mode of operation of the 3-RPS}

We study here the degenerate DKP for the first operation mode $\mathfrak{I}_1$ of the 3-RPS (where rotation part is a half-turn). The analysis for the second mode of operation would be analogous.

The boundary $\mathfrak{J}_1$ of the mode of operation has been described in section \ref{bord3RPS} as a projective 2-plane with homogeneous coordinates, $\left(w_1, w_2, w_3\right)$. The direction of translation is given by:
\[
\mathbf{u}=\begin{bmatrix}
w_1^2+w_2^2+w_3^2\\
0\\
0
 \end{bmatrix}\]
 and we normalize the homogeneous coordinates so that $w_1^2+w_2^2+w_3^2=1$.
The degenerate DKP is to solve the system
\[
\begin{aligned}
d_1&=3k_2w_1w_3+\sqrt{3} k_2w_1w_2\\
d_2&=2\sqrt{3} k_2w_1w_2\\
w_1^2+w_2^2+w_3^2&=1\;.
\end{aligned}
\]
for $w_1,w_2,w_3$; actually opposite solutions to the system give the same solution of the degenerate DKP.
\begin{center}
\includegraphics[scale=.3]{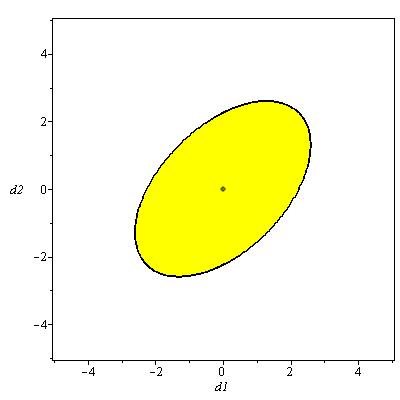}
\end{center}
The picture above illustrates that there are two solutions to the degenerate DKP inside the ellipse (whose equation is $d_1^2+d_2^2-d_1d_2=81/16$) except for the origin, and no solution outside the ellipse. The set of critical values is the ellipse plus the origin. At the origin, the degenerate DKP has infinitely many real solutions, actually one solution with $w_2=w_3=0$ (half-turn with vertical axis) and a projective line of solutions $w_1=0$ (all half-turns with horizontal axes). The latter projective line is the intersection of $\mathfrak{J}_1$ with the other boundary $\mathfrak{J}_2$, and this intersection of boundaries of modes of operation is indeed a degenerate self-motion of the platform.

We compare now the degenerate DKP with the DKP for larger and larger fixed values of $r_3$.  We have used the Maple library SIROPA to represent the number of real solutions in function of $r_1,r_2$, and also the critical values of the IKM. The pictures below represent the computed set of critical values of the IKM and the number of real solutions to the DKP for $r_3=5,10,50$. These pictures of slices can be compared with Figures 5 and 8 in \cite{SWH14}.
\begin{center}
\includegraphics[scale=.25]{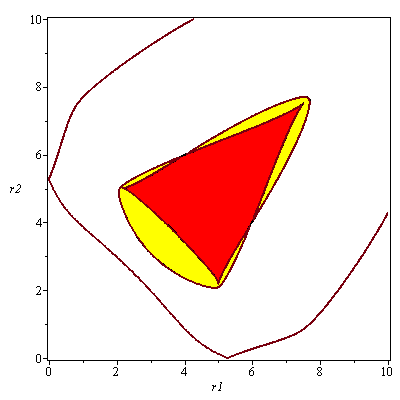}
\includegraphics[scale=.25]{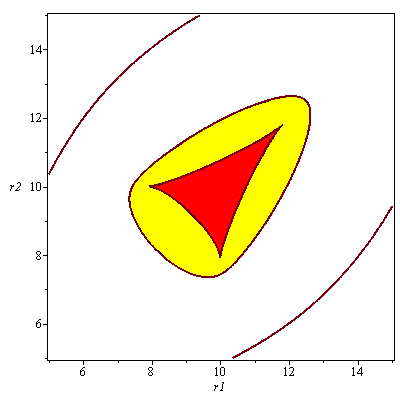}
\includegraphics[scale=.25]{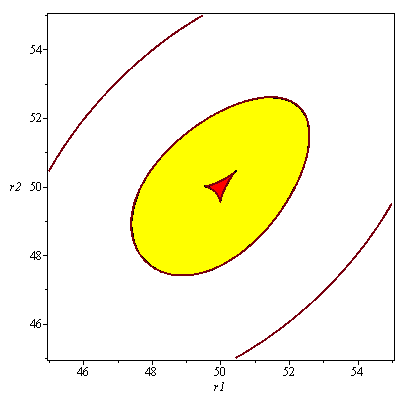}
\end{center}
There are 4 solutions to the DKP in the yellow region (inside the  small oval and outside the deltoid) and 8 in the red one (inside the deltoid). The deltoid is shrinking to a point as $r_3$ tends to infinity, giving in the limit the picture of the degenerate DKP. There is also an oval of the curves of critical values which is not relevant for the analysis of the number of real solutions to the DKP; it disappears for the degenerate DKP. 

It remains to explain the discrepancy between the number of solutions (4 for the DKP and 2 for the degenerate DKP in the yellow region). Take for instance the solutions for $d_1=2, d_2=1$ in the degenerate DKP:
$$\begin{aligned}
&w_1=0.90.., w_2=0.21..,w_3=0.36..\\
&w_1=0.42.., w_2=0.45..,w_3=0.78..
\end{aligned}
$$
and the solutions for $r_1=52, r_2=51,r_3=50$ in the DKP: 
$$\begin{aligned}
&x_1=0.90.., x_2=0.20.., x_3=0.37.., z=50.9..\\
&x_1=-0.90.., x_2=0.20.., x_3=0.37.., z=-50.9..\\
&x_1=0.42.., x_2=0.45.., x_3=0.78.., z=50.9..\\
&x_1=-0.42.., x_2=0.45.., x_3=0.78.., z=-50.9..
\end{aligned}
$$
where $z$ is the height of the center of the platform. In the solutions to the DKP, there are two pairs of configurations symmetric with respect to the base plane. But since the vector $\mathbf{u}$ is always pointing upwards, the solutions to the DKP where the platform is under the base plane give in the limit opposite $d_1$ and $d_2$ in the formulas \ref{eq:degIKM}.

\subsection{Degenerate DKP for the mysterious mode of operation of the Tsai 3-UPU}

We now turn to the mode of operation $\mathfrak{K}_8$ of the Tsai 3-UPU. Very few is known about the kinematic analysis of this mode. It is explained in \cite{WH01} that the DKP for this mode is of degree 64, and that there are values for $r_i$ with 24 real solutions. We show how to obtain some more information concerning the DKP by degenerating it. 

We have seen in section \ref{sec: bdryconf} that the boundary of this mode decomposes into $\mathfrak L_{5}\cup \mathfrak L_{6}\cup \mathfrak L_{7}$. We study separately the degenerate DKP for each of these components.

For $\mathfrak{L}_5$, we use homogeneous coordinates $(w_1,w_2,w_3)$ that we normalize with $w_1^2+w_2^2+w_3^2=1$, we take $y_0=1$ and all other coordinates $w_0,y_1,y_2,y_3$ are 0. The rotation part is the half-turn with axis directed by $(w_1,w_2,w_3)$ and the direction of translation is given by 
$$\mathbf{u}= \begin{bmatrix}
w_1\\w_2\\ w_3
\end{bmatrix}\;.$$
The equations for the degenerate DKP are then simply:
\[
\begin{aligned}
4d_1&=\sqrt3 w_2 -3w_3 \\
2d_2&=\sqrt{3} w_2\\
w_1^2+w_2^2+w_3^2&=1.
\end{aligned}
\]
 
For $\mathfrak{L}_6$, we use homogeneous coordinates $(w_0,w_2,w_3)$ that we normalize with $w_0^2+w_2^2+w_3^2=1$, we take $y_1=1$ and all other coordinates $y_0,w_1,y_2,y_3$ are 0. The rotation part is a rotation with horizontal axis, and the direction of translation is given by 
$$
\mathbf{u}=\begin{bmatrix}
-w_0\\
-w_3\\
w_2
 \end{bmatrix}\;.$$
The equations for the degenerate DKP are:
\[
\begin{aligned}
4d_1&=-15 w_2 +5\sqrt3 w_3 \\
2d_2&=5\sqrt{3} w_3\\
w_0^2+w_2^2+w_3^2&=1.
\end{aligned}
\]
Finally, $\mathfrak{L}_7$ has only two real points for which $d_1=d_2=0$. Hence, we obtain the following picture for the degenerate DKP:

\begin{center}
\includegraphics[scale=.5]{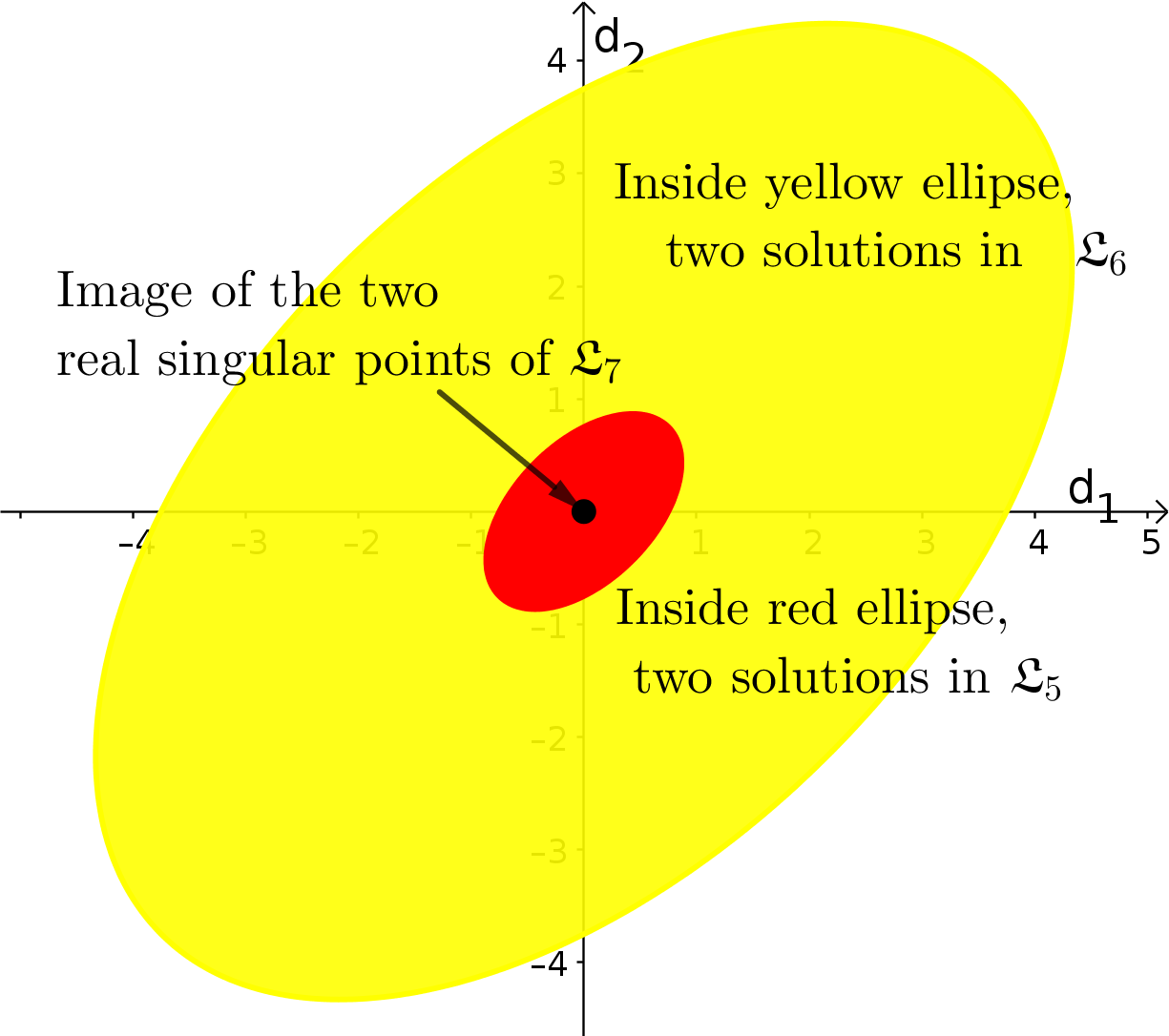}
\end{center}

One has to take a little care in the interpretation of this picture. Indeed, the two systems of coordinates $(w_1,w_2,w_3)$ and $(-w_1,-w_2,-w_3)$ represent the same point in $\mathfrak{L}_5$, but the corresponding $(d_1,d_2)$ have opposite signs. The same applies to $\mathfrak{L}_6$.

The comparison with the DKP of this mode of operation with its degenerate version can only be performed partially due to the complexity of the DKP. We compare  what happens for $d_2=0$ and $d_1$ variable and what happens for $r_2=r_3=100$ and $r_1$ variable. For the degenerate DKP, the (rather heavy) computation made with SIROPA shows that there are 0 solution for $|d_1|>15/4$, 2 solutions for $3/4<|d_1|<15/4$ and 4 solutions for $0<|d_1|<3/4$. For the DKP, there are 0 solution for $r_1-100<-3.75..$ or $r_1-100>3.74..$, 2 solutions for $-3.75..<r_1-100<-0.77..$ or $0.72..<r_1-100<3.74..$, 4 solutions for $-0.77..<r_1-100<-0.18..$ or $0.48..<r_1-100<0.72..$; for $r_1$ closer to $100$, there are up to 12 solutions.\par \medskip
The separation inside the mode $\mathfrak{K}_8$ between the two degenerate modes $\mathfrak{L}_5$ and $\mathfrak{L}_6$ can be well seen for $r_1=100.6, r_2=r_3=100$. The rotation part of the four solutions to the DKP are given by $x_0 = \pm0.98.., x_1 = 0., x_2 = 0.15.., x_3 = 0.$ and $x_0 = 0., x_1 = \pm0.73.., x_2 = 0., x_3 = 0.68..$.

In conclusion, the degenerate DKP is very easy and provides some relevant information on the very difficult DKP, for sufficiently large lengths of limbs.

\bibliographystyle{amsplain} 
\bibliography{Bibliocompact}

\end{document}